\documentclass{amsart}
\usepackage{amsfonts}
\usepackage{amsmath}
\usepackage{enumerate}
\usepackage{amsthm}
\usepackage{hyperref}
\usepackage{cite}

\begin{document}

\newtheorem{definition}{Definition}[section]
\newtheorem{theorem}[definition]{Theorem}
\newtheorem{proposition}[definition]{Proposition}
\newtheorem{remark}[definition]{Remark}
\newtheorem{lemma}[definition]{Lemma}
\newtheorem{corollary}[definition]{Corollary}
\newtheorem{example}[definition]{Example}

\numberwithin{equation}{section}

\title[Prescribed Chern scalar curvature flow...]{Prescribed Chern scalar curvature flow on compact Hermitian manifolds with negative Gauduchon degree}
\author{Weike Yu}


\begin{abstract}
In this paper, we present a unified flow approach to prescribed Chern scalar curvature problem on compact Hermitian manifolds with negative Gauduchon degree. When the conformal class of its Hermitian metric contains a balanced metric, we give some sufficient conditions on the candidate curvature function $f$ which guaranties the convergence of the flow to a conformal Hermitian metric whose Chern scalar curvature is $f$.
\end{abstract}
\maketitle

\section{introduction}
Let $(M^n, J, g)$ be a compact Hermitian manifold with complex dimension $n\geq 2$, where $J$ is a complex structure and $g$ is a Hermitian metric on $M$. The fundamental form $\omega$ of $g$ is given by $\omega(\cdot,\cdot)=g(J\cdot, \cdot)$, and we will confuse the Hermitian metric $g$ and its fundamental form $\omega$ in this paper. On a Hermitian manifold $(M, \omega)$, there is a unique linear connection $\nabla^{Ch}$ (called the Chern connection) preserving the Hermitian metric and the complex structure, whose torsion has vanishing $(1,1)$-part everywhere. The scalar curvature with respect to $\nabla^{Ch}$, which is referred to as the Chern scalar curvature, can be given by
\begin{align}
S^{Ch}(\omega)=\text{tr}_{\omega}\sqrt{-1}\bar{\partial}\partial\log\omega^n,
\end{align}
where $\omega^n$ is the volume form of $\omega$. If $\tilde{\omega}=e^{\frac{2}{n}u}\omega$ with $u\in C^\infty(M)$ is a conformal Hermitian metric to $\omega$, according to \cite{[ACS]}, the Chern scalar curvature of $\tilde{\omega}$ is given by 
\begin{align}\label{1.2..}
S^{Ch}(\tilde{\omega})=e^{-\frac{2}{n}u}\left(-\Delta^{Ch}_{\omega}u+S^{Ch}(\omega)\right),
\end{align}
where $\Delta_{\omega}^{Ch}$ denotes the Chern Laplacian with respect to $\omega$. 

The prescribed Chern scalar curvature problem is a Hermitian analogue of prescribing scalar curvatures in Riemannian geometry: Given a smooth real-valued function $f$ on $M$, does there exist a Hermitian metric $\tilde{\omega}$ conformal to $\omega$ (i.e., $\tilde{\omega}=e^{\frac{2}{n}u}\omega$) such that its Chern scalar curvature $S^{Ch}(\tilde{\omega})=f$? According to \eqref{1.2..}, it is equivalent to solving the following equation:
\begin{align}
-\Delta^{Ch}_{\omega}u+S^{Ch}(\omega)=fe^{\frac{2}{n}u}. 
\end{align}
When $f$ is a constant, the above question is referred to as the Chern-Yamabe problem, which is was first proposed by Angella-Calamai-Spotti in \cite{[ACS]}. Later, the prescribed Chern scalar curvature problem was investigated by several authors, see \cite{[CZ], [DL], [Fus], [Ho], [HS], [LM], [Ma], [SZ], [WZ2], [Yu1], [Yu2], [Yuan]} and so on.

In this paper, we will give a unified flow approach to the prescribed Chern scalar curvature problem on the compact Hermitian manifold $(M, \omega_0)$ with negative Gauduchon degree, which includes two cases: $f\leq (\not\equiv0)$ and $f$ is sign-changing. Inspired by \cite{[RR]}, we consider the following prescribed Chern scalar curvature flow on $(M, \omega_0)$:
\begin{equation}
\left\{
\begin{aligned}
&\partial_t \omega(t)=-\left(S^{Ch}(\omega(t))-f\right)w(t)\\
&\omega(0)=w_0,
\end{aligned}
\right.
\label{1.1}
\end{equation}
where the initial data $\omega_0=e^{\frac{2}{n}u_0}\eta$ for some smooth function $u_0\in C^\infty(M)$, $\eta$ is the unique Gauduchon metric in the conformal class $\{\omega_0\}$ with volume $1$ (cf. Theorem \ref{theorem2.2}) satisfying the Gauduchon degree (see \eqref{2.6.} for its definition)
\begin{align}\label{1.2.}
\Gamma(\{\omega_0\})=\int_MS^{Ch}(\eta)d\mu_{\eta}<0.
\end{align}
Since equation \eqref{1.1} preserves the conformal structure of $M$, then any smooth solution of \eqref{1.1} is of the form $\omega(t)=e^{\frac{2}{n}u(t)}\eta$, where $u(t)\in C^\infty(M)$. In terms of $u(t)$, \eqref{1.1} can be written in the following equivalent form:

\begin{equation}
\left\{
\begin{aligned}
&\partial_t e^{\frac{2}{n}u}=\Delta^{Ch}_{\eta}u-S_0+fe^{\frac{2}{n}u}\\
&u(0)=u_0,
\end{aligned}
\right.
\label{1.2}
\end{equation}
where $S_0$ denotes the Chern scalar curvature of the Gauduchon metric $\eta$. 

Our first theorem gives the global existence of the flow \eqref{1.2} for any $f\in C^\infty(M)$ on the compact Hermitian manifold with negative Gauduchon degree.

\begin{theorem}\label{thm1.1}
Let $(M^n, \omega_0)$ be a compact Hermitian manifold with Gauduchon degree $\Gamma(\{\omega_0\})<0$. Then for any function $f\in C^\infty(M)$ and any $u_0\in C^\infty(M)$, the flow \eqref{1.2} has a unique global smooth solution $u\in C^\infty(M\times [0,\infty))$.
\end{theorem}

Our main result of this paper is the following theorem:

\begin{theorem}\label{thm1.2}
Let $(M^n, \omega_0)$ be a compact Hermitian manifold with Gauduchon degree $\Gamma(\{\omega_0\})<0$ and the conformal class $\{\omega_0\}$ containing a balanced metric. Suppose that there exists a function $u^*\in C^2(M)$ such that 
\begin{align}\label{1.4}
-\Delta^{Ch}_\eta u^*+S_0-fe^{\frac{2}{n}u^*}\geq 0.
\end{align}
Then for any initial data $u_0\in C^\infty(M)$ with $u_0\leq u^*$, the solution of \eqref{1.2} converges to $u_\infty$ in $C^\infty(M)$ as $t\rightarrow +\infty$, where $u_\infty$ satisfies
\begin{align}
-\Delta^{Ch}_\eta u_\infty+s_0=fe^{\frac{2}{n}u_\infty}.
\end{align}
\end{theorem}

Since in Section \ref{sec5} we prove that when the candidate function $f\leq 0\ (\not\equiv 0)$, there exists a $C^2(M)$ function $u^*$ satisfying \eqref{1.4}, then we have the following corollary:
\begin{corollary}
Let $(M^n, \omega_0)$ be a compact Hermitian manifold with Gauduchon degree $\Gamma(\{\omega_0\})<0$ and the conformal class $\{\omega_0\}$ containing a balanced metric. Then for any smooth function $f\leq 0\ (\not\equiv0)$ and any initial data $u_0\in C^\infty(M)$, the solution of \eqref{1.2} converges to $u_\infty$ in $C^\infty(M)$ as $t\rightarrow +\infty$, where $u_\infty$ satisfies
\begin{align}
-\Delta^{Ch}_\eta u_\infty+s_0=fe^{\frac{2}{n}u_\infty}.
\end{align}
\end{corollary}
Note that the initial data $u_0$ in the above corollary can be taken as an arbitrary smooth function on $M$, because the parameters $b$ in the super-solution $u^*$ we constructed in \eqref{5.9} can be arbitrarily large.

A direct consequence of the above corollary gives
\begin{corollary}
Let $(M^n, \omega_0)$ be a compact Hermitian manifold with Gauduchon degree $\Gamma(\{\omega_0\})<0$ and the conformal class $\{\omega_0\}$ containing a balanced metric. Then for any smooth function $f\leq 0\ (\not\equiv 0)$, there exists a conformal Hermitian metric $\hat{\omega}\in\{\omega_0\}$ such that its Chern scalar curvature $S^{Ch}(\hat{\omega})=f$.
\end{corollary}

Next we consider the case that the candidate function $f$ is sign-changing. Since $M$ is compact, every sign-changing function $f\in C^\infty(M)$ can be expressed as $f=f_0+\lambda$, where $f_0$ is a nonconstant smooth function satisfying $\max_M f_0=0$ and $\lambda>0$ is a constant. In section \ref{sec5}, we show that when $\lambda$ is a small positive number, there exists a $C^2(M)$ function $u^*$ satisfying \eqref{1.4} for such $f=f_0+\lambda$. Therefore, we obtain the following corollary, which provides a flow approach to Theorem 1.5 in our recent paper \cite{[Yu1]}.

\begin{corollary}
Let $(M^n, \omega_0)$ be a compact Hermitian manifold with Gauduchon degree $\Gamma(\{\omega_0\})<0$ and the conformal class $\{\omega_0\}$ containing a balanced metric. Then for any smooth function $f=f_0+\lambda$, where $f_0$ is a nonconstant smooth function satisfying $\max_M f_0=0$ and $\lambda>0$ is a small number whose upper bound is given by \eqref{5.14..}, there exists a conformal Hermitian metric $\hat{\omega}\in\{\omega_0\}$ such that its Chern scalar curvature $S^{Ch}(\hat{\omega})=f$.
\end{corollary}
\begin{remark}
When $(M^n, \omega_0)$ is a compact Riemannian surface (i.e., $\dim_{\mathbb{C}}M=1$), the authors of \cite{[BEW]} gave a different flow approach to the prescribed Gauss curvature problem with sign-changing candidate function $f=f_0+\lambda$.
\end{remark}

This paper is organized as follows. In Section 2, we recall some basic notions and notations related to the prescribed Chern scalar curvature problem. In Section 3, we prove the global existence of the flow. In Section 4, we study the long-time behavior of the flow and prove our main result Theorem \ref{thm1.2}. In Section 5, we will construct super-solution of prescribed Chern scalar curvature equation under some sufficient conditions on the candidate curvature function $f$.

\textbf{Acknowledgements.} The author would like to thank Prof. Yuxin Dong and Prof. Xi Zhang for their continued support and encouragement.

\section{Preliminaries}

Suppose that $(M^n, J, g)$ is a Hermitian manifold with complex dimension $n$ and its fundamental form $\omega(\cdot,\cdot)=g(J\cdot, \cdot)$. Let $TM^{\mathbb{C}}=TM\otimes \mathbb{C}$ be the complexified tangent space of $M$, and we extend $J$ and $g$ from $TM$ to $TM^{\mathbb{C}}$ by $\mathbb{C}$-linearity. Then we have the decomposition:
\begin{align}
TM^{\mathbb{C}}=T^{1,0}M\oplus T^{0,1}M,
\end{align}
where $T^{1,0}M$ and $T^{0,1}M$ are the eigenspaces of complex structure $J$ corresponding to the eigenvalues $\sqrt{-1}$ and $-\sqrt{-1}$, respectively. Note that every $m$-form can also be decomposed into $(p, q)$-forms for each $p, q \geq 0$ with $p + q = m$ by extending the complex structure $J$ to forms.

On a Hermitian manifold $(M^n, J, g)$, there exists a unique affine connection $\nabla^{Ch}$ preserving both the Hermitian metric $g$ and the complex structure $J$, that is,
\begin{align}
 \nabla^{Ch}g=0,\quad \nabla^{Ch}J=0, 
\end{align}
whose torsion $T^{Ch}(X,Y)=\nabla_XY-\nabla_YX-[X,Y]$ satisfies 
\begin{align}
T^{Ch}(JX,Y)=T^{Ch}(X,JY)
\end{align}
for any $X, Y\in TM$. In this paper, we will confuse the Hermitian metric $g$ and its corresponding fundamental form $\omega$. 

For a Hermitian manifold $(M^n, \omega)$ with Chern connection $\nabla^{Ch}$, it is well-known that the Chern scalar curvature of $\omega$ is given by
\begin{align}
S^{Ch}(\omega)=\text{tr}_{\omega}Ric^{(1)}(\omega)=\text{tr}_{\omega}\sqrt{-1}\bar{\partial}\partial \log{\omega^n},
\end{align}
where $Ric^{(1)}(\omega)$ is the first Chern Ricci curvature.

In Hermitian geometry, there is a canonical elliptic operator called the Chern Laplace operator: for any smooth function $u: M\rightarrow \mathbb{R}$, the Chern Laplacian $\Delta^{Ch}$ of $u$ is defined by
\begin{align}
\Delta^{Ch}u=-2\sqrt{-1} tr_{\omega}\overline{\partial}\partial u.
\end{align}
\begin{lemma}[cf. \cite{[Gau]}]\label{lemma2.1}
On a compact Hermitian manifold $(M^n,\omega)$, we have
\begin{align}
-\Delta^{Ch}u=-\Delta_d u+(du,\theta)_{\omega},
\end{align}
where $\Delta_d$ is the Hodge-de Rham Laplacian, $\theta$ is the Lee form or torsion $1$-form given by $d\omega^{n-1}=\theta\wedge \omega^{n-1}$, and $(\cdot, \cdot)_\omega$ denotes the inner product on $1$-form induced by $\omega$. Furthermore, if $\omega$ is balanced, namely $\theta\equiv0$, then $-\Delta^{Ch}u=-\Delta_d u$ for any $u\in C^\infty(M)$.
\end{lemma}
Let 
\begin{align}
\{\omega\}=\{e^{\frac{2}{n}u}\omega\ |\ u\in C^\infty(M)\}
\end{align}
denote the conformal class of the Hermitian metric $\omega$. In \cite{[Gau1]}, Gauduchon proved the following theorem:
\begin{theorem}\label{theorem2.2}
Let $(M, \omega)$ be a compact Hermitian manifold with $\dim_{\mathbb{C}} M\geq 2$, then there exists a unique Gauduchon metric $\eta\in \{\omega\}$ (i.e., $d^*\theta=0$) with volume $1$.
\end{theorem}
In terms of the above theorem, one can define an invariant $\Gamma(\{\omega\})$ of the conformal class $\{\omega\}$ which is called the Gauduchon degree:
\begin{align}\label{2.6.}
\Gamma(\{\omega\})=\frac{1}{(n-1)!}\int_M c^{BC}_1(K^{-1}_M)\wedge\eta^{n-1}=\int_M S^{Ch}(\eta)d\mu_\eta,
\end{align}
where $\eta$ is the unique Gauduchon metric in $\{\omega\}$ with volume $1$, $c^{BC}_1(K^{-1}_M)$ is the first Bott-Chern class of anti-canonical line bundle $K^{-1}_M$, and $d\mu_\eta$ denotes the volume form of the Gauduchon metric $\eta$.

Fix a compact Hermitian manifold $(M^n, \omega)$, we consider the conformal change $\widetilde{\omega}=e^{\frac{2}{n}u}\omega$. From \cite{[Gau]}, the Chern scalar curvatures of $\widetilde{\omega}$ and $\omega$ have the following relationship:
\begin{align}\label{2.7}
-\Delta^{Ch}_\omega u+S^{Ch}(\omega)=S^{Ch}(\widetilde{\omega})e^{\frac{2}{n}u},
\end{align}
where $S^{Ch}(\omega)$ and $S^{Ch}(\widetilde{\omega})$ denote the Chern scalar curvatures of $\omega$ and $\widetilde\omega$, respectively. 

To resolve the prescribed Chern scalar curvature problem, it is useful to see the equation \eqref{2.7} as an Euler-Lagrange equation of some functional. Unfortunately, in general, such a functional dose not always exist. In fact, according to \cite[Prop. 5.3]{[ACS]} and \cite[Prop. 2.12]{[Fus]}, we know the following
\begin{proposition}
Let $(M, \omega)$ be a compact Hermitian manifold. Then \eqref{2.7} can be seen as an Euler-Lagrange equation for standard $L^2$ pairing if and only if $\omega$ is balanced.
\end{proposition}

\section{Global existence of the flow}
In this section, we will prove the global existence of the flow \eqref{1.2}. Since the flow \eqref{1.2} is a parabolic equation, then there exists a smooth solution of \eqref{1.2} defined on a maximal interval $[0, T_{\text{max}})$ with $0<T_{\text{max}}<\infty$. 
\begin{lemma}
Let $u(x,t)$ be a smooth solution of \eqref{1.2} defined on $M\times [0,T_{\text{max}})$. Then we have 
\begin{align}\label{3.7}
u(x,t)\geq -C_0=\min\left\{\min_M (u_0-v_0), \frac{n}{2}\log \frac{-\overline{S_0}}{\|f\|_{L^\infty(M)}e^{\frac{2}{n}\max_M v_0}}\right\}+\min_M v_0
\end{align}
and 
\begin{align}\label{3.1.}
u(x,t)\leq \max\{0, \max_M u_0\}+C_1t+\max_M v_0-\min_Mv_0.
\end{align}
for any $(x,t)\in M\times [0,T_{\text{max}})$, where 
\begin{align}
C_1=\frac{n}{2}\left(\|f\|_{L^\infty(M)}-\overline{S_0}\right)
\end{align}
and $v_0\in C^\infty(M)$ be a solution of 
\begin{align}
\Delta^{Ch}_\eta v_0=S_0-\overline{S_0},
\end{align}
here $\overline{S_0}=\int_M S_0d\mu_\eta<0$.
\end{lemma}

\proof Let $w=u-v_0$. From \eqref{1.2}, it follows that
\begin{align}\label{3.3}
e^{\frac{2}{n}v_0}\partial_t e^{\frac{2}{n}w}=\Delta^{Ch}_{\eta}w-\overline{S_0}+fe^{\frac{2}{n}v_0}\cdot e^{\frac{2}{n}w}.
\end{align}
For \eqref{3.7}, let $(x_0, t_0)\in M\times [0, T]$ such that $w(x_0, t_0)=\min_{M\times [0, T]}w$ for any $T\in (0,T_{\text{max}})$. If $t_0=0$, then 
\begin{align}
u(x,t)=w(x,t)+v_0(x)\geq  \min_{M\times [0, T]}w+\min_M v_0=\min_{M}(u_0-v_0)+\min_M v_0,
\end{align}
thus \eqref{3.7} is proved in this case. If $t_0>0$, then by the maximum principle, we have $\partial_t w(x_0, t_0)\leq 0$ and $\Delta^{Ch}_{\eta}w(x_0,t_0)\geq 0$. Thus, by \eqref{3.3},
\begin{equation}
\begin{aligned}
0&\geq -\overline{S_0}-\|f\|_{L^\infty(M)}e^{\frac{2}{n}\max_M v_0}e^{\frac{2}{n}w(x_0,t_0)},
\end{aligned}
\end{equation}
which implies that
\begin{equation}
\begin{aligned}
u(x,t)&=w(x,t)+v_0(x)\\
&\geq  \min_{M\times [0, T]}w+\min_M v_0\\
&\geq \frac{n}{2}\log \frac{-\overline{S_0}}{\|f\|_{L^\infty(M)}e^{\frac{2}{n}\max_M v_0}}+\min_M v_0.
\end{aligned}
\end{equation}
Therefore, \eqref{3.7} is proved. For the inequality \eqref{3.1.}, we set $v=w-C_1t$, where $C_1=\frac{n}{2}\left(\|f\|_{L^\infty(M)}-\overline{S_0}\right)$. Let $(x_1, t_1)\in [0, T]$ such that $v(x_1, t_1)=\max_{M\times [0, T]}v$ for any $T\in (0,T)$. If $t_1=0$, then 
\begin{equation}
\begin{aligned}
u(x,t)&=w(x,t)+v_0(x)\\
&=v(x,t)+C_1t+v_0(x)\\
&\leq \max_{M\times [0, T]}v+C_1t+ v_0(x)\\
&=(u_0(x_1)-v_0(x_1))+C_1t+v_0(x)\\
&\leq \max_M u_0-\min_Mv_0+C_1t+\max_M v_0.
\end{aligned}
\end{equation}
If $t_1>0$, by the maximum principle we have $\partial_t v(x_1, t_1)\geq 0$ and $\Delta^{Ch}_{\eta}v(x_1,t_1)\leq 0$, i.e., $\partial_t w(x_1, t_1)\geq C_1$ and $\Delta^{Ch}_{\eta}w(x_1,t_1)\leq 0$. Substituting these in \eqref{3.3} yields
\begin{align}
\frac{2}{n}C_1e^{\frac{2}{n}u(x_1,t_1)}\leq -\overline{S_0}+\|f\|_{L^\infty(M)}e^{\frac{2}{n}u(x_1,t_1)}
\end{align}
which implies that
\begin{align}
e^{\frac{2}{n}u(x_1,t_1)}\leq 1,
\end{align}
that is, $u(x_1,t_1)\leq 0$, because of $C_1=\frac{n}{2}\left(\|f\|_{L^\infty(M)}-\overline{S_0}\right)$. Therefore,
\begin{equation}
\begin{aligned}
u(x,t)&=w(x,t)+v_0(x)\\
&=v(x,t)+C_1t+v_0(x)\\
&\leq v(x_1,t_1)+C_1t+v_0(x)\\
&=w(x_1,t_1)+C_1(t-t_1)+v_0(x)\\
&=u(x_1,t_1)-v_0(x_1)+C_1(t-t_1)+v_0(x)\\
&\leq-\min_M v_0+ C_1t+\max_Mv_0.
\end{aligned}
\end{equation}
which finishes the proof of \eqref{3.1.}.
\qed

\proof[Proof of Theorem \ref{thm1.1}] Let $u$ be the solution of \eqref{1.2} defined on a maximal interval $[0, T_{\text{max}})$. Assume by contradiction that $T_{\text{max}}<+\infty$. Applying the Krylov-Safonov estimate \cite{[KS]} to \eqref{1.2}, we get the H\"order estimate: there exists $\delta\in(0,1)$ and $C>0$ independent of $t\in[0,T]$ such that
\begin{equation*}
\begin{aligned}
\|&u\|_{C^{\frac{\delta}{2}, \delta}(M\times [0, T])}\\
&\leq C\left(\|u\|_{L^\infty(M\times [0, T])}+\left\|-\frac{n}{2}S_0e^{-\frac{2}{n}u}+f\right\|_{L^\infty(M\times [0, T])}\right)\\
&\leq C_3(u_0, f, C_0, S_0, T),
\end{aligned}
\end{equation*}
where $C_3(u_0, f, C_0, S_0, T)$ is positive constant depending only on $u_0, f, C_0, S_0, T$. Applying the classical theory of parabolic equations, we get $\|u\|_{C^k(M\times [0,T])}\leq C_k$ for any $k\in\mathbb{N}$, where $C_k>0$ is a constant depending only on $u_0, f, C_0, S_0, T$. Hence, we can extend the solution $u$ beyond $T_{\text{max}}$, which contradictes the maximality of $T_{\text{max}}$. Therefore, $T_{\text{max}}=+\infty$, which finishes the proof of this theorem.
\qed

\section{Long-time behavior of the flow}
In this section, we will study the long-time behavior of the flow \eqref{1.2} and prove our main result Theorem \ref{thm1.2}. Using the method of \cite{[RR]}, we prove that

\begin{lemma}\label{lem4.1}
Let $(M, \omega_0)$ be a compact Hermitian manifold with Gauduchon degree $\Gamma(\{\omega_0\})<0$. Suppose that there exists a function $u^*\in C^2(M)$ such that 
\begin{align}\label{4.1}
-\Delta^{Ch}_\eta u^*+S_0-fe^{\frac{2}{n}u^*}\geq 0.
\end{align}
Then for any initial data $u_0\in C^\infty(M)$ with $u_0\leq u^*$, the solution $u(x,t): M\times [0,+\infty)\rightarrow \mathbb{R}$ of \eqref{1.2} satisfies 
\begin{align}
u(x,t)\leq u^*(x)
\end{align}
for any $(x,t)\in M\times [0,+\infty)$. Combining with \eqref{3.7}, we have the uniform $C^0$-estimate:
\begin{align}\label{4.3}
\|u\|_{C^0(M\times [0, +\infty))}\leq C_5,
\end{align}
where $C_5>0$ is a constant only dependent on $M,\omega_0, S_0, \|f\|_{L^\infty(M)}, u^*$.
\end{lemma}
\proof 
Let $v(x,t)=u^*(x)-u(x,t)$. From \eqref{1.2} and \eqref{4.1}, it follows that
\begin{align}\label{4.4}
\partial_t\left(e^{\frac{2}{n}u^*}-e^{\frac{2}{n}u} \right)\geq \Delta^{Ch}_{\eta}v+f\left(e^{\frac{2}{n}u^*}-e^{\frac{2}{n}u}\right).
\end{align}
Set 
\begin{align}
a(x,t)=\frac{2}{n}\int_0^1e^{\frac{2}{n}\{su^*(x)+(1-s)u(x,t)\}}ds>0,
\end{align}
then we have
\begin{align}
av=e^{\frac{2}{n}u^*}-e^{\frac{2}{n}u}.
\end{align}
Thus, it follows from \eqref{4.4} that
\begin{align}
\partial_t\left(av \right)\geq \Delta^{Ch}_{\eta}v+fav,
\end{align}
which gives
\begin{align}
\partial_t v \geq a^{-1}\Delta^{Ch}_{\eta}v+\left(f-a^{-1}\partial_ta\right)v.
\end{align}
In terms of the maximum principle \cite[Lemma 2.3, page 8]{[Lieb]} and $v(\cdot, 0)\geq 0$, we get $v(x,t)\geq 0$ for any $t\in [0,+\infty)$, i.e., $u(x,t)\leq u^*(x)$. 
\qed

Now we suppose that the conformal class $\{\omega_0\}$ contains a balanced metric, then the Gauduchon metric $\eta$ is balanced. Define the functional
\begin{align}\label{4.9}
E(u)=\frac{1}{2}\int_M|\nabla u|^2_\eta d\mu_\eta+\int_MS_0ud\mu_\eta-\frac{n}{2}\int_Mfe^{\frac{2}{n}u}d\mu_\eta,
\end{align}
then we have
\begin{lemma}
Let $(M^n, \omega_0)$ be a compact Hermitian manifold with Gauduchon degree $\Gamma(\{\omega_0\})<0$ and the conformal class $\{\omega_0\}$ containing a balanced metric. Then the functional $E$ is non-increasing along the flow \eqref{1.2}, i.e.,
\begin{equation}\label{4.10}
\begin{aligned}
\frac{d}{dt}E(u(\cdot,t))=-\frac{2}{n}\int_M| \partial_tu|^2e^{\frac{2}{n}u}d\mu_\eta\leq 0.
\end{aligned}
\end{equation}
\end{lemma}
\proof Since $\eta$ is balanced, a direct computation gives
\begin{equation}
\begin{aligned}
\frac{d}{dt}E(u(\cdot,t))&=-\int_M\Delta^{Ch}_\eta u\cdot \partial_tud\mu_\eta+\int_MS_0\partial_tud\mu_\eta-\int_Mfe^{\frac{2}{n}u}\partial_tud\mu_\eta\\
&=\int_M\left(-\Delta^{Ch}_\eta u+S_0-fe^{\frac{2}{n}u}\right)\partial_tud\mu_\eta\\
&=-\frac{2}{n}\int_M| \partial_tu|^2e^{\frac{2}{n}u}d\mu_\eta\leq 0.
\end{aligned}
\end{equation}
\qed

\proof[Proof of Theorem \ref{thm1.2}]

Applying the Krylov-Safonov estimate \cite{[KS]} to \eqref{1.2} and using Lemma \ref{lem4.1}, we obtain the uniform H\"older estimates for the solution of the flow: there exists $\delta\in(0,1)$ and $C>0$ independent of $t\in[0, +\infty)$ such that

\begin{equation*}
\begin{aligned}
\|&u\|_{C^{\frac{\delta}{2}, \delta}(M\times [0, +\infty))}\\
&\leq C\left(\|u\|_{L^\infty(M\times [0, +\infty))}+\left\|-\frac{n}{2}S_0e^{-\frac{2}{n}u}+f\right\|_{L^\infty(M\times [0, +\infty))}\right)\\
&\leq C_7,
\end{aligned}
\end{equation*}
where $C_7$ is a positive constant independent of $t\in[0,+\infty)$. According to the classical regularity theory for linear parabolic equations, for any $k\in\mathbb{N}$, we obtain 
\begin{align}\label{4.12}
\|u\|_{C^k(M\times [0,+\infty))}\leq C_k,
\end{align}
where $C_k$ is a positive constant independent of $t\in[0,+\infty)$. On the other hand, from \eqref{4.10}, it follows that
\begin{align}
\int^s_0\int_M| \partial_tu|^2e^{\frac{2}{n}u}d\mu_\eta dt=-\frac{n}{2}\int^s_0\frac{d}{dt}E(u(\cdot,t))dt=\frac{n}{2}\left(E(u_0)-E(u(\cdot, t)) \right).
\end{align}
Next, we claim that there exists a constant $C_8>0$ independent of $t\in [0,+\infty)$ such that
\begin{align}
E(u(\cdot, t))\geq -C_8
\end{align}
for any $t\in [0,+\infty)$. Indeed, by \eqref{3.7}, \eqref{4.3} and \eqref{4.9}, we get
\begin{equation}
\begin{aligned}
E(u(\cdot, t))&\geq \int_MS_0(u+C_0)d\mu_\eta -C_0\int_MS_0d\mu_\eta -\frac{n}{2}\|f\|_{C^0(M)}e^{\frac{2}{n}C_5}\\
&\geq \inf_MS_0\int_M(u+C_0)d\mu_\eta -C_0\int_MS_0d\mu_\eta -\frac{n}{2}\|f\|_{C^0(M)}e^{\frac{2}{n}C_5}\\
&\geq \inf_MS_0C_5+ \inf_MS_0C_0-C_0\int_MS_0d\mu_\eta -\frac{n}{2}\|f\|_{C^0(M)}e^{\frac{2}{n}C_5}\\
&>-\infty,
\end{aligned}
\end{equation}
where $\inf_MS_0<0$ since $\Gamma(\{\omega_0\})=\int_MS_0d\mu_\eta<0$. Therefore, 
\begin{align}
\int^{+\infty}_0\int_M| \partial_tu|^2e^{\frac{2}{n}u}d\mu_\eta dt<+\infty.
\end{align}
Hence, there is a sequence $t_k\rightarrow +\infty$ such that
\begin{align}\label{4.17}
\lim_{k\rightarrow +\infty}\int_M| \partial_tu(\cdot, t_k)|^2e^{\frac{2}{n}u(\cdot, t_k)}d\mu_\eta=0.
\end{align}
Since $M$ is compact, it follows from \eqref{4.12} that by passing to a subsequence (if necessary), $u(\cdot, t_k)$ converges to a smooth function $u_\infty$ in $C^m(M)$ for any $m\in\mathbb{N}$. Replacing $t$ by $t_k$ in \eqref{1.2} and letting $k\rightarrow +\infty$, from \eqref{4.3} and \eqref{4.17},  we deduce that $u_\infty$ satisfies
\begin{align}
-\Delta^{Ch}_\eta u_\infty+s_0=fe^{\frac{2}{n}u_\infty}.
\end{align}
As in \cite{[RR]}, applying the general result of Simon \cite[Corollary 2]{[Sim]}, we obtain that 
\begin{align}
\lim_{t\rightarrow+\infty}\|u(\cdot, t)-u_\infty\|_{C^k(M)}=0
\end{align}
for any $k\in\mathbb{N}$, which finishes the proof of this theorem.
\qed

\section{The construction of super-solutions}\label{sec5}
In this section, we will provide some sufficient conditions on $f$ which can guarantee the existence of the super-solution of 
\begin{align}\label{5.1..}
-\Delta^{Ch}_\eta u+S_0=fe^{\frac{2}{n}u}
\end{align}
on the compact Hermitian manifold $(M, \omega_0)$ with Gauduchon degree $\Gamma(\{\omega_0\})<0$, where $\eta$ is the unique Gauduchon metric in $\{\omega_0\}$ with volume $1$. By Integrating on M, it follows from \eqref{5.1..} that
\begin{align}
\int_Mfe^{\frac{2}{n}u}d\mu_\eta=-\int_M\Delta^{Ch}_\eta ud\mu_\eta+\int_MS_0d\mu_\eta=\Gamma(\{\omega_0\})<0,
\end{align}
because $\eta$ is a Gauduchon metric, which implies the necessary condition:
\begin{align}
\min_Mf<0.
\end{align}
Hence, we consider the following three cases.\\
\textbf{Case 1}: $f\leq 0\ (\not\equiv 0)$.

Let $v_0\in C^\infty(M)$ be a solution of 
\begin{align}\label{5.2}
\Delta^{Ch}_\eta v_0=S_0-\overline{S_0},
\end{align}
where $\overline{S_0}=\int_MS_0d\mu_\eta=\Gamma(\{\omega_0\})<0$, and let $v_1\in C^\infty(M)$ be a solution of
\begin{align}
\Delta^{Ch}_\eta v_1=\overline{f}-f,
\end{align}
where $\overline{f}=\int_M fd\mu_\eta<0$, since $f\leq 0\ (\not\equiv 0)$. Set
\begin{align}
u^*=v_0+av_1+b,
\end{align}
where $a>0$ and $b\in\mathbb{R}$ are two constants that we will choose later, and let 
\begin{align}
c^-_0=\inf_Mv_0, \quad c^-_1=\inf_M v_1,
\end{align}
by \eqref{5.2}, we obtain
\begin{align}
-\Delta^{Ch}_\eta u^*+S_0-f e^{\frac{2}{n} u^*}\geq\overline{S_0}-a\overline{f}+f\left(a-e^{\frac{2}{n}(c^-_0+ac^-_1+b)}\right)\geq0,
\end{align}
provided that $a, b$ satisfy
\begin{equation}\label{5.9}
\left\{
\begin{aligned}
&a\geq\frac{\overline{S_0}}{\overline{f}},\\
&b\geq\frac{n}{2}\ln a-c^-_0-ac^-_1.\\
\end{aligned}
\right.
\end{equation}

Next we consider the case that the candidate curvature function $f$ is sign-changing. Since $M$ is compact, then every sign-changing function $f\in C^\infty(M)$ can be expressed as $f=f_0+\lambda$, where $f_0$ is a nonconstant smooth function satisfying $\max_M f_0=0$ and $\lambda>0$ is a constant. Hence, we consider the following case:\\
\textbf{Case 2}: $f=f_0+\lambda$, where $\max_Mf_0=0$ and $f_0\not\equiv0$, and $\lambda>0$ is a small constant.

Let $v_2\in C^\infty(M)$ be a solution of
\begin{align}\label{5.10.}
\Delta^{Ch}_\eta v_2=\overline{f_0}-f_0, \quad v_2>0.
\end{align}
where $\overline{f_0}=\int_M f_0d\mu_\eta<0$, since $\max_Mf_0=0$ and $f_0\not\equiv0$. Set
\begin{align}
u^*=v_0+av_2+b,
\end{align}
where $a>0, b\in\mathbb{R}$ are two constants that we will choose later, and let 
\begin{align}
c^+_0=\sup_Mv_0,\quad c_2^-=\inf_M v_2>0,\quad c^+_2=\sup_Mv_2,
\end{align}
then by \eqref{5.2} and \eqref{5.10.}, we deduce that
\begin{equation*}
\begin{aligned}
-\Delta^{Ch}_\eta u^*+S_0-f e^{\frac{2}{n} u^*}&=\overline{S_0}-a\overline{f_0}-\lambda e^{\frac{2}{n}(v_0+av_2+b)}+f_0\left(a-e^{\frac{2}{n}(v_0+av_2+b)}\right)\\
&\geq \overline{S_0}-a\overline{f_0}-\lambda e^{\frac{2}{n}(c^+_0+ac^+_2+b)}+f_0\left(a-e^{\frac{2}{n}(c^-_0+ac^-_2+b)}\right)\\
&\geq 0,
\end{aligned}
\end{equation*}
provided that 
\begin{equation}
\left\{
\begin{aligned}
&\overline{S_0}-a\overline{f_0}> 0,\\
&a-e^{\frac{2}{n}(c^-_0+ac^-_2+b)}\leq 0,\\
&\overline{S_0}-a\overline{f_0}-\lambda e^{\frac{2}{n}(c^+_0+ac^+_2+b)}\geq 0,\\
\end{aligned}
\right.
\end{equation}
which is equivalent to
\begin{equation}\label{5.14..}
\left\{
\begin{aligned}
&a>\frac{\overline{S_0}}{\overline{f_0}},\\
&e^{\frac{2}{n}b}\geq ae^{-\frac{2}{n}(c^-_0+ac^-_2)},\\
&0<\lambda\leq(\overline{S_0}-a\overline{f_0})e^{-\frac{2}{n}(c^+_0+ac^+_2+b)}. \\
\end{aligned}
\right.
\end{equation}
\textbf{Case 3}: Suppose that $\dim_{\mathbb{C}}M=1$ and $f\in C^\infty(M)$ satisfies the following condition:
\begin{align}\label{5.15.}
\frac{1}{\|fe^{\frac{2}{n}v_0}\|_{L^\infty(M)}}\int_M f^+e^{\frac{2}{n}v_0}d\mu_\eta\leq C_M\left(\frac{1}{\|fe^{\frac{2}{n}v_0}\|_{L^\infty(M)}}\int_M f^-e^{\frac{2}{n}v_0}d\mu_\eta\right)^\theta
\end{align}
where $f^+(x)=\max\{f(x),0\}$, $f^-=\max\{-f(x),0\}$, $C_M>0$ is constant depending only on $M$, $\theta=\frac{\pi-2\pi\chi(M)+1}{\pi-1}$, $\chi(M)$ is the Euler characteristic of $M$ and $v_0$ is defined by \eqref{5.2}.

To give a super-solution $u^*$ of \eqref{5.1..}, it is sufficient to construct a $C^2(M)$ function $U^*$ such that
\begin{align}\label{5.6}
-\Delta^{Ch}_\eta U^*+\overline{S_0}-ge^{\frac{2}{n}U^*}\geq 0,
\end{align}
where $g=fe^{\frac{2}{n}v_0}$. Indeed, Let $u^*=U^*+v_0$, by \eqref{5.2} and \eqref{5.6}, we obtain
\begin{align}
-\Delta^{Ch}_\eta u^*+S_0-fe^{\frac{2}{n}u^*}=-\Delta^{Ch}_\eta U^*+\overline{S_0}-ge^{\frac{2}{n}U^*}\geq 0.
\end{align}
By \eqref{5.15.}, $g$ satisfies 
\begin{align}
\frac{1}{\|g\|_{L^\infty(M)}}\int_M g^+d\mu_\eta\leq C_M\left(\frac{1}{\|g\|_{L^\infty(M)}}\int_M g^-d\mu_\eta\right)^\theta.
\end{align}
According to \cite[Theorem 2]{[RR]}, the above condition leads to the existence of $U^*$ mentioned above, and thus we get the super-solution $u^*$. Therefore, we can recover the all results of \cite{[RR]} when $\dim_{\mathbb{C}}M=1$.

\bigskip
\bigskip

Weike Yu

School of Mathematical Sciences, 

Ministry of Education Key Laboratory of NSLSCS,

Nanjing Normal University,

Nanjing, 210023, Jiangsu, P. R. China,

wkyu2018@outlook.com

\bigskip

\end{document}